\documentclass[11pt,a4paper]{article}
\usepackage[latin1]{inputenc}
\usepackage{graphicx}
\usepackage[mathscr]{euscript}
\usepackage{a4}
\usepackage{epsfig}
\usepackage{stmaryrd,amsfonts,amsmath,amssymb}
\usepackage{ifthen}

\newcommand{\origforall}{}    \let\origforall=\forall
\renewcommand{\forall}{\,\origforall\,}
\newcommand{\origexists}{}    \let\origexists=\exists
\renewcommand{\exists}{\,\origexists\,}
\newcommand{\origin}{}        \let\origin=\in
\renewcommand{\in}{{\,\origin\,}}
\newcommand{\origwedge}{}     \let\origwedge=\wedge
\renewcommand{\wedge}{{\quad\origwedge\quad}}
\newcommand{\origvee}{}       \let\origvee=\vee
\renewcommand{\vee}{{\quad\origvee\quad}}
\newcommand{\origfrac}{}      \let\origfrac=\frac
\renewcommand{\frac}[2]{{\;\origfrac{#1}{#2}\;}}
    
\renewcommand{\subset}{\subseteq}

\newcommand{\beq}{\begin{eqnarray*}}
\newcommand{\eeq}{\end{eqnarray*}}

\newcommand{\K}{\mathbb{K}}

\newcommand{\N}{\mathbb{N}}
\newcommand{\Z}{\mathbb{Z}}

\newcommand{\B}{\mathfrak{B}}
\newcommand{\n}{\textnormal}

\DeclareMathOperator{\ord}{ord}
\DeclareMathOperator{\nord}{nord}
\DeclareMathOperator{\Inn}{Inn}
\DeclareMathOperator{\Env}{Env}
\newcommand{\id}{\mathrm{id}}
\DeclareMathOperator{\ev}{ev}
\DeclareMathOperator{\Alg}{Alg}
\DeclareMathOperator{\End}{End}
\newcommand{\my}{\mu}
\newcommand{\trid}{\triangleright}
\newcommand{\sq}[2]{Q_{{#1},{#2}}}  
\newcommand{\fg}[1]{C_{#1}}

\newcommand{\Perm}[1]{S_{#1}} 
\newcommand{\APerm}[1]{A_{#1}} 
\newcommand{\op}{\mathrm{op}}

\newcommand{\na}[1]{ %
\ifthenelse{\equal{#1}{3A}}{\B(Q_{3,1},-1)}{ %
\ifthenelse{\equal{#1}{3B}}{\B(Q_{3,1},E_3)^{(2)}}{ %
\ifthenelse{\equal{#1}{4A}}{\B(Q_{4,1},-1)^{(2)}}{ %
\ifthenelse{\equal{#1}{4B}}{\B(Q_{4,1},-1)^{(\neq 2)}}{ %
\ifthenelse{\equal{#1}{4C}}{\B(Q_{4,1},\chi_4)}{ %
\ifthenelse{\equal{#1}{5A}}{\B(Q_{5,2},-1)}{ %
\ifthenelse{\equal{#1}{5B}}{\B(Q_{5,3},-1)}{ %
\ifthenelse{\equal{#1}{6A}}{\B(Q_{6,1},-1)}{ %
\ifthenelse{\equal{#1}{6B}}{\B(Q_{6,1},\chi_6)}{ %
\ifthenelse{\equal{#1}{6C}}{\B(Q_{6,2},-1)}{ %
\ifthenelse{\equal{#1}{7A}}{\B(Q_{7,4},-1)}{ %
\ifthenelse{\equal{#1}{7B}}{\B(Q_{7,5},-1)}{ %
\ifthenelse{\equal{#1}{10A}}{\B(Q_{10,1},-1)}{ %
\ifthenelse{\equal{#1}{10B}}{\B(Q_{10,1},\chi_{10})}{ %
\ifthenelse{\equal{#1}{15A}}{\B(Q_{15,7},-1)}{ %
\ifthenelse{\equal{#1}{5*}}{\B(Q_{5,*},-1)}{ %
\ifthenelse{\equal{#1}{6*}}{\B(Q_{6,*},*)}{ %
\ifthenelse{\equal{#1}{7*}}{\B(Q_{7,*},-1)}{ %
\ifthenelse{\equal{#1}{10*}}{\B(Q_{10,1},*)}{ %
\textnormal{UNKNOWN NICHOLS ALGEBRA}}}}}}}}}}}}}}}}}}}}}

\newcommand{\esquare}{\hfill\ensuremath{\square}}

\newcounter{claim}
\newcounter{titleclaim}
\newenvironment{Proof}{{\textbf{Proof}$\;$}}{\esquare\newline}
\newenvironment{Non-Proof}{{\textbf{Non-Proof}$\;$}}{\hfill\ensuremath{\times}\newline}

\newtheorem{LemmaBase}[claim]{Lemma}
\newenvironment{Lemma}[1][-1]{\ifthenelse{\equal{#1}{-1}}
    {\begin{LemmaBase}}{\begin{LemmaBase}[#1]}
    \dotfill \, {\em \bf \arabic{claim}} \quad\\{} }
    {\addtocounter{titleclaim}{1}\end{LemmaBase}}

\newtheorem{CorollaryBase}[claim]{Corollary}
\newenvironment{Corollary}[1][-1]{\ifthenelse{\equal{#1}{-1}}
    {\begin{CorollaryBase}}{\begin{CorollaryBase}[#1]}
    \dotfill \, {\em \bf \arabic{claim}} \quad\\{} }
    {\addtocounter{titleclaim}{1}\end{CorollaryBase}}

\newtheorem{PropositionBase}[claim]{Proposition}
\newenvironment{Proposition}[1][-1]{\ifthenelse{\equal{#1}{-1}}
    {\begin{PropositionBase}}{\begin{PropositionBase}[#1]}
    \dotfill \, {\em \bf \arabic{claim}} \quad\\{} }
    {\addtocounter{titleclaim}{1}\end{PropositionBase}}

\newtheorem{ConjectureBase}[claim]{Conjecture}

\newtheorem{ExampleBase}[claim]{Example}
\newenvironment{Example}[1][-1]{\ifthenelse{\equal{#1}{-1}}
    {\begin{ExampleBase}}{\begin{ExampleBase}[#1]}
    \dotfill \, {\em \bf \arabic{claim}} \quad\\{} }
    {\addtocounter{titleclaim}{1}\end{ExampleBase}}

\newtheorem{DefinitionBase}[claim]{Definition}
\newenvironment{Definition}[1][-1]{\ifthenelse{\equal{#1}{-1}}
    {\begin{DefinitionBase}}{\begin{DefinitionBase}[#1]}
    \dotfill \, {\em \bf \arabic{claim}} \quad\\{} }
    {\addtocounter{titleclaim}{1}\end{DefinitionBase}}

\newtheorem{NotationBase}[claim]{Notation}

\newtheorem{TheoremBase}[claim]{Theorem}
\newenvironment{Theorem}[1][-1]{\ifthenelse{\equal{#1}{-1}}
    {\begin{TheoremBase}}{\begin{TheoremBase}[#1]}
    \hrulefill \, {\em \bf \arabic{claim}} \quad\\{} }
    {\addtocounter{titleclaim}{1}\end{TheoremBase}}

\newtheorem{DLemmaBase}[claim]{Dependend Lemma}

\newtheorem{DCorollaryBase}[claim]{Dependend Corollary}

\newtheorem{DPropositionBase}[claim]{Dependend Proposition}

\newtheorem{DTheoremBase}[claim]{Dependend Theorem}

\begin{document}

\title{Divisibility Relations for the Dimensions and Hilbert series of
Nichols Algebras of Non-Abelian Group Type}
\author{Andreas Lochmann}
\maketitle

\abstract
We present a divisibility relation for the dimensions and Hilbert series
of certain classes of Nichols algebras of non-abelian group type, which
generalizes Nichols algebras over Coxeter groups with constant cocycle $-1$. For this we
introduce three groups of isomorphisms acting on Nichols algebras,
which generalizes the exchange operator introduced by Milinski and
Schneider in \cite{MiS} for Coxeter groups.

\section{Introduction}

In \cite{GHV_quadratic}, Theorem 4.14, Gra\~na, Heckenberger and Vendramin
gave a full classification of finite-dimensional Nichols algebras of
non-abelian group type with absolutely irreducible Yetter-Drinfeld modules
under the assumption that their Hilbert series factorize as
\begin{equation*}
\mathcal{H}_\B(t) \;\;=\;\; \prod_{j=1}^r \, (\alpha_j)_t
\quad \textnormal{with} \quad (\alpha_j)_t \,:=\,
1\,+\,t\,+\,t^2\,+\,\cdots\,+\,t^{\alpha_j-1}
\end{equation*}
for some $r,\,\alpha_j\,\in\,\N$.
From the examples in Subsection \ref{SEC_tables}, they covered all Nichols
algebras except $\na{3B}$ and $\na{4C}$, whose Hilbert series do not factorize as above.
In a subsequent paper (\cite{HLV_cubic}), Heckenberger, Vendramin and the
author classified all finite-dimensional Nichols algebras of non-abelian
group type with absolutely irreducible Yetter-Drinfeld modules under the
more general assumption 
\begin{equation} \label{EQ_cubic_hilbert}
\mathcal{H}_\B(t) \;\;=\;\; \prod_{j=1}^r \, (\alpha_j)_t \,\cdot\,
\prod_{j=1}^s \, (\beta_j)_{t^2} \,, \qquad \alpha_j,\,\beta_j\,\in\,\N\,,
\end{equation}
but restricting the calculations to the special class of braided racks.
With this approach, the Nichols algebras $\na{3B}$ and $\na{4C}$ were found
to be finite-dimensional.
Though the calculations are intricate already in the case of braided racks,
the classification of finite-dimensional Nichols algebras 
along the approach proposed by Gra\~na, Heckenberger and Vendramin
seems to be feasible; if it is possible to show that the Hilbert series always
factorizes in a way similar to Equation \ref{EQ_cubic_hilbert}.
This paper wants to contribute to the last question.

Shortly after publication of \cite{HLV_cubic}, Heckenberger brought to 
our attention that the dimensions of the Nichols algebras identified
so far are always divisible by the order of the inner group of their
underlying racks, a fact which has been proven for Nichols algebras
over Coxeter groups by Milinski and Schneider in 1999 (\cite{MiS}).
For their proof, they constructed a family of isomorphisms between the
homogeneous components, based on interchanging the coefficients
of a decomposition of each vector. This decomposition is possible
due to the Nichols-Zoeller-Theorem (\cite{NZ_grouplike}, see also
Theorem 7.2.9 in \cite{DNR}), which essentially states that a finite-dimensional
Hopf algebra $H$ is a free left $B$-module for all Hopf subalgebras $B$
of $H$.

Similar freeness theorems are abundant in Hopf theory and can be found
e.g.\ in \cite{AHS_semisimple} (Theorem 1 and Lemma 3.2),
\cite{Skr_freeness} (Theorem 6.1), and \cite{Rad_subalgebras}.
We will use a version by Gra\~na:

\begin{Theorem}[Gra\~na] \label{THE_grana_freeness}
(Cf. \cite{G_freeness}, Theorem 3.8.1) Let $V$ be a finite-dimensional
Yetter-Drinfeld module over a group $G$ and assume $V\,=\,V'\,\oplus\,U$ with
\begin{itemize}
\item
$V'\,\subset\,V$ a $G'$-stable $\K G$-subcomodule, where $G'$ is the smallest
subgroup of $G$ with $\delta(V')\,\subseteq\,\K G'\,\otimes\,V'$, and
\item
$U\,\subset\,V$ a $\K G$-subcomodule and $\K G'$-submodule.
\end{itemize}
Let $\{e_i\,:\,i\,=\,1,\ldots, d\}$ be a basis of $V'$ and $\partial_i$
the corresponding braided (skew) derivatives on the Nichols algebra $\B(V)$ over $V$. 
Then
$$ \B(V)\;\;\cong\;\;\left(\bigcap_{i=1}^d\,\ker\partial_i\right)\,\otimes\,\B(V') $$
as right $\B(V')$-modules and left $(\bigcap_{i=1}^d\,\ker\partial_i)$-modules.

\noindent
In particular, $\dim\B(V')$ divides $\dim\B(V)$.
\end{Theorem}

Applying this together with specialized maps similar to those used by
Milinski and Schneider, our main results are as follows:

\begin{Theorem} \label{THM_PRE_inner_group_and_subalgebra}
Let $\B$ be a finite-dimensional Nichols algebra over the indecomposable quandle $X$ with a
2-cocycle $\chi$. Assume that the degree of $X$ divides the order of the diagonal elements of $\chi$.
Then each $\Inn X$-homogeneous component of $\B$ has the same dimension
(this is wrong for $\deg X\,\nmid\,\ord \chi$ in general) and thus $\#\Inn X$
divides $\dim \B$.

Moreover, let $X'$ be a non-empty proper sub-rack of $X$ and $\B'$
the Nichols sub-algebra generated by $X'$. Assume that $X\setminus X'$
still generates $\Inn X$. Then $\#\Inn X\cdot \dim\B'$ divides $\dim\B$.
\end{Theorem}

If we drop the assumption that the degree of $X$ divides the order of the diagonal
elements of $\chi$, we can still prove the following:

\begin{Theorem} \label{THM_PRE_hilbert}
Let $\B$ be a finite-dimensional Nichols algebra over an indecomposable
rack $X$ and a 2-cocycle with diagonal elements of order $m$. Let $X'$ be a non-empty proper
subrack of $X$ and $\B'$ its corresponding Nichols sub-algebra of $\B$.
Then the Hilbert series $\mathcal{H}_\B(t)$ is divisible by
$(m)_t\cdot \mathcal{H}_{\B'}(t)$.
\end{Theorem}

Our paper is organized as follows. In Section \ref{SEC_preliminaries}, we
define the core notions to understand the main results (racks, quandles,
Nichols algebras, (opposite) braided derivations, Hilbert series) as well as a
list of known Nichols algebras of non-abelian group type with absolutely
irreducible Yetter-Drinfeld modules and the corresponding quandles
in Subsection \ref{SEC_tables}. In Section \ref{SEC_shifts}, we introduce
maps, so-called shifts, which are direct generalisations of the maps
defined by Milinski and Schneider in \cite{MiS}.
An even more generalized form 
is currently considered by Angiono, Vay, and Vendramin.
In Subsection \ref{SEC_n_divides_m},
we apply them to show the first half of Theorem \ref{THM_PRE_inner_group_and_subalgebra},
i.e.~$\#\Inn X\,|\,\dim\B$ for a slightly larger class
of Nichols algebras than those examined in \cite{MiS}.
In Section \ref{SEC_modified_shifts}, we define two improved variations
of shifts and apply them to Gra\~na's Freeness Theorem to show
the second half of 
Theorem \ref{THM_PRE_inner_group_and_subalgebra} for the same class
of Nichols algebras as before. We will then use our methods to analyze
arbitrary Nichols algebras over non-trivial, indecomposable racks and
proof Theorem \ref{THM_PRE_hilbert}.
Finally, we will show in Subsection \ref{SEC_nichols_4B}, that the
direct approach to proof $\#\Inn X\,|\,\dim\B$ along the lines of \cite{MiS}
and Subsection \ref{SEC_n_divides_m} is not feasible in the general
case $\deg X\,\nmid\,\ord \chi$. 

\section{Preliminaries} \label{SEC_preliminaries}

\subsection{Notations}

Denote $\N\,=\,\N_0\setminus\{0\}$. With $C_k$ we denote the cyclic
group or order $k$, and
with $[j]_k$ the equivalence class of $j$ in $C_k$, for $k\,\in\,\N$ and
$j\in\Z$. We use $\mathfrak{S}_n$ to refer to the symmetric group on
$n$ symbols.
Furthermore, we use the notation ``$\sq{x}{y}$'' to refer to the $y$-th
indecomposable quandle of size $x$ in Gra\~na's and Vendramin's list of
small indecomposable quandles (\cite{V_racks}, implemented in Rig, see \cite{Rig}).
$\K$ shall always be a field of arbitrary characteristic, if not said
otherwise.

\subsection{Nichols Algebras from Racks}

For a detailled account on racks in the context of Nichols Algebras,
see \cite{AG_racks}.

\begin{Definition} \label{DEF_rack}
A {\em rack} $X$ is a set with a binary operation $\trid$, which fulfills:
\begin{itemize}
\item
Left-self-distributivity: For all $t_1,t_2,t_3\in X$ holds
$t_1\,\trid\,(t_2\,\trid\, t_3)\,=\,(t_1\,\trid\, t_2)\,\trid\,(t_1\,\trid\, t_3)$.
\item
The operations $g_t\,:\,X\,\rightarrow\,X$, $s\,\mapsto\,t\,\trid\, s$ are bijections.
\end{itemize}
An idempotent rack is called {\em quandle}.

Due to left-self-distributivity, each $g_t$ as defined above is an automorphism of $(X,\,\trid)$.
The permutation subgroup generated by the $g_t$ is called the {\em inner group} $\Inn X$ of $X$.
It is a quotient of the enveloping (or structure) group
$$ \Env X\;:=\;\langle t\,\in\,X \;|\; s\cdot t\,=\,(s\,\trid t)\cdot s
\; \forall\, s,\,t\,\in\,X\rangle_\n{group}\,. $$
A rack $X$ is called {\em indecomposable}, if $\Inn X$ acts transitively on $X$.
If the map $g_\cdot:\,X\,\rightarrow\,\Inn X$, $t\,\mapsto\,g_t$ is injective,
$X$ is called {\em faithful}.
\end{Definition}

If $X$ is a faithful quandle, then $X$ is realized as a conjugation-closed
generating subset of a group. On the other hand,
each such subset is a faithful quandle.

Throughout this article, let $X$ denote a finite indecomposable faithful quandle
and $\K$ our ground-field.

\begin{Example}
(Cf. \cite{AG_racks}, subsection 1.3)
Let $A$ be an abelian group and $\alpha:\,A\,\rightarrow\,A$ some automorphism of $A$.
Then $\trid:\,A\times A\,\rightarrow\,A$, $(t,\,s)\,\mapsto\,\alpha(t\,-\,s)\,+\,s$
defines a quandle structure on $A$. Quandles of this kind
are called {\em affine quandles}. Many of the quandles used to construct Nichols
algebras are affine (see the tables in subsection \ref{SEC_tables}), though not all.
Assume $A$ is an affine quandle with two commuting elements $t,\,s\,\in\,A$
(i.e.\ $t\,\trid\,s\,=\,s$ and vice versa), then inserting into the definition
gives $t\,=\,s$. This excludes many quandles, e.g.\ the quandles
given by transpositions in the symmetric group $S_n$ for $n\,\geq\,4$.
The smallest non-affine quandles without commuting elements are $Q_{15,5}$
and $Q_{15,6}$ (\cite{Rig}). On the other hand, there are many non-faithful
affine quandles.
\end{Example}

\begin{Definition} \label{DEF_nichols_algebra}
(Cf.~\cite{MiS}; for terms of Hopf algebra theory, see e.g.\ \cite{DNR})
Let $V_0$ be a finite-dimensional vector space and set $V\,=\,\K X\,\otimes\,V_0$.
Denote $V_t\,:=\,\{t\}\,\otimes\,V_0$.
Assume the inner group $\Inn X$ acts on $V$ with
$g_t\,(V_s)\,\subseteq\,V_{t\trid s}$ for all $t,\,s\,\in\,X$. Then $V$ is
an $\Inn X$-Yetter-Drinfeld module, and a braiding
$$ c(v\,\otimes\,w) \;\;=\;\; g_t(w)\,\otimes\,v $$
on the tensor product $V\,\otimes\,V$ is induced,
where $v\,\in\,V_t$ and $w\,\in\,V$ are arbitrary.
This in turn induces a co-algebra structure on the tensor algebra $(TV,\,\mu)$,
which is uniquely determined by the two requirements
\begin{enumerate}
\item
$\Delta\,\mu\;\;=\;\;(\mu\,\otimes\,\mu)\,(\id_V\,\otimes\,c\,\otimes\,\id_V)\,(\Delta\,\otimes\,\Delta)$
\item
and each $v\,\in\,V\,\subset\,TV$ is primitive.
\end{enumerate}
The co-unit is given by $\epsilon(1)\,:=\,1$ and $\epsilon(v)\,:=\,0$ for all $v\,\in\,V\,\subset TV$.
Furthermore, an antipode $S:\,TV\,\rightarrow\,TV$ can be defined which
endows $TV$ with the structure of an $\N_0$-graded Hopf algebra
(actually $F_X$-graded, where $F_X$ is the free group over the set $X$).

There is a unique maximal homogeneous ideal and coideal $I$ of $TV$ such that
$$ \Delta(I)\;\;\subseteq\;\;I\,\otimes\,TV\,+\,TV\,\otimes\,I $$
and such that all homogeneous elements of $I$ are of degree $\geq\,2$. The quotient
$\B\,:=\,TV/I$ is called the {\em Nichols algebra} of $V$.
It is an $\Env X$-graded braided Hopf algebra, whose primitive elements
are exactly those of degree $1$ and generate $\B$.
\end{Definition}

The Nichols algebra $\B$ can be completely described in terms of the rack
$X$ and a 2-cocycle $\chi:\,X\,\times\, X\,\rightarrow\,\End(V_0)$ satisfying
$$ \chi(t,\,s\trid r)\,\chi(s,\,r)\;\;=\;\;\chi(t\trid s,\,t\trid r)\,\chi(t,\,r)\,, $$
which induces the action of $\Inn X$ on $V$ (\cite{AG_racks}, \cite{AFGV_simple}).

If $V$ is finite-dimensional and absolutely irreducible as a Yetter-Drinfeld module,
the Lemma of Schur shows that $\chi(t,t)$ actually is a scalar multiple $q_t$
of the identity for each $t\,\in\,X$ (\cite{GHV_quadratic}; Theorem 2.7 in \cite{G_low}).
If $X$ is indecomposable, the transitive action of $\Inn X$ ensures that
$q_t$ does not depend on $t$; we drop the index in this case.
However, Gra\~na showed in \cite{G_low}, Lemma 3.1, that the cases $\dim V_0\,\geq\,2$
impose severe restrictions on $q$ and $\chi$ if $\B$ is to be finite-dimensional. 
Therefore, for the most part of this paper, we will restrict to the case
$\dim V_0\,=\,1$, without losing too many cases.

\begin{Definition}
Let $X$ be an indecomposable rack. Then the order of $g_t$ does not depend on
$t\,\in\,X$ and we define the degree $\deg X$ to be this number.

\noindent
The nilpotency order $\nord v$ of an element $v\in\B$ is the minimal $m\,\in\,\N$ with $v^m\,=\,0$.

\noindent
The order $\ord \chi$ of a 2-cocycle $\chi$ is the minimal $m\,\in\,\N\,\cup\,\{\infty\}$ with $\chi(t,t)^m\,=\,1$
for all $t\,\in\,X$.
(If $\B$ is finite dimensional, one easily shows that $q_t\,:=\,\chi(t,t)$ has to be a
root of unity, such that $\ord \chi$ is finite.)
\end{Definition}

For the rest of this paper, set $n\,:=\,\deg X$ and $m\,:=\,\ord \chi$.

Denote with $\{e_t\}_{t\in X}$ the standard base in $\K X$. If $V_0$ is
one-dimensional, we use it as standard base for $V$ as well; else, we
denote $e_t\otimes v$ with $e_t\,v$ (see \cite{AFGV_simple}).
If $X$ is indecomposable, $\nord (e_t\,v)$ does not depend on $t\,\in\,X$
nor $v\,\in\,V_0$.
Using Equation \ref{EQ_pld_8} in Proposition \ref{PRO_left_derivative}
one easily calculates $\nord (e_t\,v)\,=\,m$ in this case for $v\,\in\,V_0\setminus\{0\}$.
There is however no obvious relation between $m$ and $n$, as can be seen in
the examples of subsection \ref{SEC_tables}.

\begin{Definition}
Given a $G$-graded vector space $U$, denote with $U(g)$ the $g$-homogeneous
subspace of $U$. We say the grading is {\em balanced}, if $\dim U(g)$ is
finite and constant for all $g\in G$. A vector space $U$
is {\em $G$-balanced}, if $U$ is $G$-graded and the grading is balanced.
\end{Definition}

Let $U$ be a $G$-graded vector space and $H$ a quotient of $G$.
Then $U$ is $H$-graded as well.
Any Nichols algebra $\B$ is $\Env X$-graded. As there are canonical surjective
homomorphisms $\Env X\,\rightarrow\,\Z$ and $\Env X\,\rightarrow\,\Inn X$,
we will use the notation $\B(x)$ for any $x\,\in\,\Env X$, $\Z$ or
$\Inn X$ without further notice, as the latter two gradings are
induced by the first one. If the $G$-grading is balanced, then
the induced $H$-grading is balanced as well.

\begin{Definition}
If $U$ is a $\Z$-graded vector space and each $U(g)$ is finite-dimensional,
define the (formal) Hilbert series $\mathcal{H}_U$ by
$$ \mathcal{H}_U(t)\;\;:=\;\;\sum_{j\in\Z}\,\dim U(j)\,t^j\,. $$
For any $k\,\in\,\N\,\cup\,\{\infty\}$ define $(k)_t$ to be the series $\sum_{j=1}^{k-1}\,t^j$.
\end{Definition}

Finite dimensional Nichols algebras of abelian group type (i.e.\ over trivial
quandles) have been completely classified by I.~Heckenberger in \cite{H_thesis}
and \cite{H_classification}. The classification of finite dimensional Nichols
algebras of non-abelian group type, particularly over indecomposable quandles,
is advanced by several strategies. A very interesting ansatz is to identify
the set of appropriate racks, so by identifying racks of type D. This led
to the exclusion of conjugacy class racks of whole classes of groups,
notably the alternating groups $A_m$ for $m\geq 6$ (\cite{AFGV_alternating})
and many sporadic groups (\cite{AFGV_sporadic}).
An alternative is to derive inequalities on the maximal dimensions of the
lower homogeneous degrees, as has been done in \cite{GHV_quadratic} and
\cite{HLV_cubic}, and connect these to a certain factorization of $\mathcal{H}_\B$
in terms of $(k)_t$ and $(k)_{t^2}$.

\subsection{Examples for Nichols Algebras} \label{SEC_tables}

There are nine indecomposable and faithful quandles known to provide examples of
finite-dimensional Nichols algebras, with the following properties.

\bigskip

\noindent
\!\!\!\!\!\!\!\!\!\!
\begin{tabular}{|c|c|cc|l|} \hline
    &          &          &        & $g_t\,\in\,\Inn X$ for generating $t\,\in\,X$ \\
$X$ & $\deg X$ & $\Inn X$ & (Size) & (as perm.\ of $X$ in cycle notation) \\ \hline
$\sq{3}{1}$  & 2 & $\Perm{3}$ & $(6)$ & $g_1=(2,3),\; g_2=(1,3)$ \\
$\sq{4}{1}$  & 3 & $\APerm{4}$ & $(12)$ & $g_1=(2,3,4),\; g_2=(1,4,3)$ \\
$\sq{5}{2}$  & 4 & $\fg{5}\rtimes\fg{4}$ & $(20)$ & $g_1=(2,4,5,3),\; g_2=(1,4,3,5)$ \\
$\sq{5}{3}$  & 4 & $\fg{5}\rtimes\fg{4}$ & $(20)$ & $g_1=(2,3,5,4),\; g_2=(1,5,3,4)$ \\
$\sq{6}{1}$  & 2 & $\Perm{4}$ & $(24)$ & $g_1=(3,5)(4,6),\;g_2=(3,6)(4,5),$ \\
             &   &            &        & \quad$g_3=(1,5)(2,6)$ \\
$\sq{6}{2}$  & 4 & $\Perm{4}$ & $(24)$ & $g_1=(3,5,4,6),\;g_3=(1,6,2,5)$ \\
$\sq{7}{4}$  & 6 & $(\fg{7}\rtimes\fg{3})\rtimes \fg{2}$ & $(42)$ & $g_1=(2,6,5,7,3,4),$ \\
             &   &            &        & \quad$g_2=(1,4,5,3,7,6)$ \\
$\sq{7}{5}$  & 6 & $(\fg{7}\rtimes\fg{3})\rtimes \fg{2}$ & $(42)$ & $g_1=(2,4,3,7,5,6),$ \\
             &   &            &        & \quad$g_2=(1,6,7,3,5,4)$ \\
$\sq{10}{1}$ & 2 & $\Perm{5}$ & $(120)$ & $g_1=(2,7)(3,5)(4,6),$ \\
             &   &            &         & \quad$g_2=(1,7)(3,8)(4,10),$ \\
             &   &            &         & \quad$g_3=(1,5)(2,8)(4,9),$ \\
             &   &            &         & \quad$g_4=(1,6)(2,10)(3,9)$ \\ \hline
\end{tabular}

\bigskip

The quandles $\sq{3}{1}$, $\sq{5}{2}$, $\sq{5}{3}$, $\sq{7}{4}$,
and $\sq{7}{5}$ are affine quandles over the cyclic abelian groups of order $\# X$,
with $\alpha$ the multiplication with $2,\,3,\,2,\,5,$ and $3$, respectively.
$\sq{4}{1}$ also is an affine quandle over $\fg{2}\times \fg{2}$ with
$\alpha\,=\,\begin{pmatrix}1&1\\0&1\end{pmatrix}$. The quandles $\sq{3}{1}$,
$\sq{6}{1}$ and $\sq{10}{1}$ can also be defined as the conjugacy classes of
transpositions in the symmetric groups $\mathfrak{S}_n$ for $n\,=\,3,\,4,\,5$,
respectively.

The following fourteen finite-dimensional Nichols algebras over $\K$
of non-abelian group type are our basic examples, sorted by quandle and dimension.

\begin{center}
\begin{tabular}{|l|c|c|c|r|l|} \hline
$\B(X,\,c)$ & $\mathrm{char}(\K)$ & $n$ & $m$ & dimension & $\mathcal{H}_\B(t)$ \\ \hline
  $\na{3A}$ & $\ast$   & $2$ & $2$ &        12 & $(2)_t^2\,(3)_t$ \\ 
  $\na{3B}$ & $2$      & $2$ & $3$ &       432 & $(3)_t\,(4)_t\,(6)_t\,(6)_{t^2}$ \\ 
  $\na{4A}$ & $2$      & $3$ & $2$ &        36 & $(2)_t^2\,(3)_t^2$ \\
  $\na{4B}$ & $\neq 2$ & $3$ & $2$ &        72 & $(2)_t^2\,(3)_t\,(6)_t$ \\
  $\na{4C}$ & $\ast$   & $3$ & $3$ &     5,184 & $(6)_t^4\,(2)_{t^2}^2$ \\
  $\na{5A}$ & $\ast$   & $4$ & $2$ &     1,280 & $(4)_t^4\,(5)_t$ \\
  $\na{5B}$ & $\ast$   & $4$ & $2$ &     1,280 & $(4)_t^4\,(5)_t$ \\
  $\na{6A}$ & $\ast$   & $2$ & $2$ &       576 & $(2)_t^2\,(3)_t^2\,(4)_t^2$ \\
  $\na{6B}$ & $\ast$   & $2$ & $2$ &       576 & $(2)_t^2\,(3)_t^2\,(4)_t^2$ \\
  $\na{6C}$ & $\ast$   & $4$ & $2$ &       576 & $(2)_t^2\,(3)_t^2\,(4)_t^2$ \\
  $\na{7A}$ & $\ast$   & $6$ & $2$ &   326,592 & $(6)_t^6\,(7)_t$ \\
  $\na{7B}$ & $\ast$   & $6$ & $2$ &   326,592 & $(6)_t^6\,(7)_t$ \\
 $\na{10A}$ & $\ast$   & $2$ & $2$ & 8,294,400 & $(4)_t^{4}\,(5)_t^2\,(6)_t^4$ \\
 $\na{10B}$ & $\ast$   & $2$ & $2$ & 8,294,400 & $(4)_t^{4}\,(5)_t^2\,(6)_t^4$ \\ \hline
\end{tabular}
\end{center}

In all cases we have $\dim V_0\,=\,1$.
$E_N$ denotes an $N$-th root of unity and the superscripts $^{(2)}$ and $^{(\neq 2)}$
refer to the field's characteristic.
The non-constant cocycles $\chi_4$, $\chi_6$ and $\chi_{10}$ are defined as follows:
\begin{eqnarray*}
\chi_4 &:=& \left(\begin{array}{rrrr}
E_3  & -E_3 & -E_3 &  E_3 \\
-E_3 &  E_3 & -E_3 &  E_3 \\
-E_3 & -E_3 &  E_3 &  E_3 \\
 E_3 &  E_3 &  E_3 &  E_3
\end{array}\right)
\end{eqnarray*}
\begin{eqnarray*}
\chi_6 &:=& \left(\begin{array}{rrrrrr}
    -1 &  1 & -1 &  1 & -1 &  1 \\
     1 & -1 &  1 & -1 & -1 &  1 \\
     1 &  1 & -1 &  1 &  1 &  1 \\
     1 &  1 &  1 & -1 &  1 &  1 \\
     1 &  1 &  1 &  1 & -1 &  1 \\
    -1 & -1 & -1 & -1 &  1 & -1
\end{array}\right) \\
\chi_{10} &:=& \left(\begin{array}{rrrrrrrrrr}
     -1&   1&   1&   1&   1&   1&   1&   1&   1&   1 \\
     -1&  -1&   1&   1&   1&   1&  -1&   1&   1&   1 \\
     -1&  -1&  -1&   1&  -1&   1&   1&  -1&   1&   1 \\
     -1&  -1&  -1&  -1&   1&  -1&   1&   1&  -1&  -1 \\
      1&   1&   1&   1&  -1&   1&  -1&  -1&   1&   1 \\
      1&   1&   1&   1&  -1&  -1&  -1&   1&  -1&  -1 \\
      1&   1&   1&   1&   1&   1&  -1&   1&   1&   1 \\
      1&   1&   1&   1&   1&   1&   1&  -1&   1&   1 \\
      1&   1&   1&   1&   1&   1&   1&   1&  -1&   1 \\
      1&   1&   1&   1&   1&   1&   1&  -1&  -1&  -1
\end{array}\right)
\end{eqnarray*}
A more concise description of $\chi_6$ and $\chi_{10}$ in terms of transpositions 
in $\mathfrak{S}_4$ and $\mathfrak{S}_5$ is given e.g.\ in \cite{MiS} (Example 5.3)
and in \cite{V_twist}.

Note that $\na{5A}$ and $\na{5B}$ as well as $\na{7A}$ and $\na{7B}$ are dual algebras
(see Example 2.1 in \cite{AFGV_simple}),
$\na{6A}$ and $\na{6B}$ as well as $\na{10A}$ and $\na{10B}$ are twist-equivalent to
each other (\cite{V_twist}).

Also note that the factorization of $\mathcal{H}_\B(t)$ in terms of $(k)_t$ and $(k)_{t^2}$
is not unique.

\subsection{Braided Derivations and Braided Commutator}

For the rest of the paper, we assume $\dim V_0\,=\,1$ for simplicity.

\begin{Definition} \label{DEF_skew-derivations}
Given the comultiplication $\Delta:\,\B\,\rightarrow\,\B$, we
can uniquely define linear maps $\partial_t$ and $\partial^\op_t:\,\B\,\rightarrow\,\B$
for arbitrary $t\,\in\,X$ via
\begin{eqnarray*}
\Delta(v)
&=& v\,\otimes\,1\;+\;\sum_{t\in X}\,\partial_t(v)\,\otimes\,e_t\;+\;
\n{some element of $\B\,\otimes\,\bigoplus_{j=2}^\infty\,\B(j)$} \\
&=& 1\,\otimes\,v\;+\;\sum_{t\in X}\,e_t\,\otimes\,\partial^\op_t(v)
\;+\;\n{some element of $\bigoplus_{j=2}^\infty\,\B(j)\,\otimes\,\B$}\,.
\end{eqnarray*}
We call these maps {\em braided derivations} and {\em opposite
braided derivations}, respectively.
\end{Definition}

The braided derivations $\partial$ and $\partial^\op$, have been introduced
by Nichols in subsection 3.3 in \cite{Nichols} under the name
``quantum differential operators''; for an account on them, we refer to
\cite{AHS_semisimple}.

\begin{Proposition} \label{PRO_left_derivative}
The maps $\partial_t$ and $\partial_t^\op$ satisfy the following properties
for all $t,\,s\,\in\,X$ and $v,\,w\,\in\,\B$:
\begin{eqnarray}
\partial_t(1) &=& 0 \label{EQ_pld_1} \\
\partial_t(e_s) &=& \delta_{t,s} \label{EQ_pld_2} \\
\partial_t(vw) &=& v \partial_t(w) \,+\, \partial_t(v) g_t(w) \label{EQ_pld_3} \\
\partial^\op_t(1) &=& 0 \label{EQ_pld_4} \\
\partial^\op_t(e_s) &=& \delta_{t,s} \label{EQ_pld_5} \\
\partial^\op_t(vw) &=& \partial^\op_t(v)\, w \,+\, v\, \partial^\op_{g_v^{-1}(t)}(w)/\chi(v,t)
\quad \n{(if $v$ is homog.)}\quad \label{EQ_pld_6} \\
\partial_s^\op\,\partial_t &=& \partial_t\,\partial_s^\op \label{EQ_pld_7} \\
\bigcap_{t\,\in\,X}\,\ker\partial_t
&=& \bigcap_{s\,\in\,X}\,\ker\partial_s^\op \;\;=\;\; \B(0) \label{EQ_pld_8}
\end{eqnarray}
(where $\delta_{t,s}$ is the Kronecker symbol.)
\end{Proposition}
\begin{Proof}
(\ref{EQ_pld_1}), (\ref{EQ_pld_4}): Obvious.

(\ref{EQ_pld_2}), (\ref{EQ_pld_5}): Follows from primitivity of $\B(1)$.

(\ref{EQ_pld_7}): Follows from co-associativity.

(\ref{EQ_pld_8}): Follows from \cite{Nichols} and \cite{G_freeness}.

(\ref{EQ_pld_3}): Let $v,\,w\,\in\,\B$ be arbitrary. Then holds:
\begin{eqnarray*}
\Delta(vw) &=& \Delta(v)\cdot\Delta(w) \\
&=&\left(v\otimes 1\,+\,\sum_{t\in X}\,\partial_t(v)\otimes e_t \,+\,\n{higher terms}\right) \\
&& \quad \cdot\left(w\otimes 1\,+\,\sum_{t\in X}\,\partial_t(w)\otimes e_t \,+\,\n{higher terms}\right) \\
&=& (vw)\otimes 1\,+\,\sum_{t\in X}\,\big(\partial_t(v)\otimes e_t\big)\cdot\big(w\otimes 1\big) \\
&& \phantom{(vw)\otimes 1}\,+\,\sum_{t\in X}\,\big(v\otimes 1\big)\cdot \big(\partial_t(w)\otimes e_t\big)
\,+\,\n{higher terms} \\
&=& (vw)\otimes 1\,+\,\sum_{t\in X}\,\big(\partial_t(v)\, g_t(w)
\,+\,v\,\partial_t(w)\big)\otimes e_t \,+\,\n{h.t.}
\end{eqnarray*}

(\ref{EQ_pld_6}): Similar to (\ref{EQ_pld_3}) we have:
\begin{eqnarray*}
&& \Delta(vw)  \;\;=\;\; \Delta(v)\cdot\Delta(w) \\
&=&\left(1\otimes v\,+\,\sum_{t\in X}\,e_t\otimes\partial^\op_t(v) \,+\,\n{higher terms}\right) \\
&& \quad \cdot\left(1\otimes w\,+\,\sum_{s\in X}\,e_s\otimes\partial^\op_s(w) \,+\,\n{higher terms}\right) \\
&=& 1\otimes (vw)\,+\,\sum_{t\in X}\,\big(e_t\otimes\partial^\op_t(v)\big)\cdot\big(1\otimes w\big) \\
&& \phantom{1\otimes (vw)}\,+\,\sum_{s\in X}\,\big(1\otimes v\big)\cdot \big(e_s\otimes\partial^\op_s(w)\big)
\,+\,\n{higher terms} \\
&=& 1\otimes (vw)\,+\,\sum_{t\in X}\big(e_t\otimes\partial^\op_t(v)\,w\big)
\,+\,\sum_{s\in X}\big(g_v(e_s)\otimes v\,\partial^\op_s(w)\big)
\,+\,\n{h.t.}
\end{eqnarray*}
where $g_v$ is defined such that $v\,\in\,\B(g_v)$. By definition we have $g_v(e_s)\,=\,\chi(v,s)\,e_{v\trid s}$
(where $v\trid s$ is short-hand for $g_v(s)$ and the 2-cocycle $\chi$ is extended in the obvious way).
Choosing $t$ such that $e_{v\trid s}\,=\,e_t$, we conclude
\begin{eqnarray*}
\Delta(vw) &=& 1\otimes (vw)\,+\,\sum_{t\in X}\,e_t\otimes\big(\partial^\op_t(v)\,w
\,+\, v\,\partial^\op_{v\trid^{-1}t}(w)/\chi(v,t)\big)
\,+\,\n{h.t.}
\end{eqnarray*}
where $v\trid^{-1}t\,:=\,g_v^{-1}(t)$.
\end{Proof}

$\partial_t$ is a right $\sigma$-skew-derivation, as one sees from Equation
\ref{EQ_pld_3}, with the endomorphism $\sigma\,=\,g_t$.
$\partial_t^\op$ is not a right skew-derivation; so we chose the word
``opposite braided derivation'', to emphasize its kinship with
$\partial_t$. Also note that $\ker\partial_t^\op\,\neq\,\ker\partial_t$ in
general.

\begin{Definition}[\cite{AS_survey}] \label{DEF_braided_commutator}
Let $\B$ be a Nichols algebra with braiding $c$. Define the {\em braided commutator}
$$ [x,\,y]_c \;:=\; \my\circ(\id\,-\,c)(x\,\otimes\, y) $$
for all $x,\,y\,\in\,\B$.
\end{Definition}

\begin{Proposition} \label{PRO_skew-derivation_and_g_commute}
For all $t\,\in\,X$ holds:
$ \partial_t\, g_t \,=\, q_t\cdot g_t\,\partial_t $.
\end{Proposition}
\begin{Proof}
By induction over the $\N_0$-degree $d$ of $v\,\in\,\B$. For $d\,\in\,\{0,\,1\}$, this
is clear. For each $v,\,w\,\,\in\,\B$ we have
\begin{eqnarray*}
\partial_t\,g_t\,(v\,w) &=& (g_t\,v)\cdot (\partial_t\,g_t\,w) \,+\, (\partial_t
\,g_t\,v)\cdot (g_t^2\, w) \\
&\stackrel{\n{ind.}}{=}& q_t\cdot\big((g_t\,v)\cdot (g_t\,\partial_t\,w) \,+\, (g_t\,\partial_t\,
v)\cdot (g_t^2\,w)\big)\\
\n{and}\quad
g_t\,\partial_t\,(v\,w) &=& (g_t\,v)\cdot (g_t\,\partial_t\,w) \,+\, (g_t\,\partial_t\,
v)\cdot (g_t^2\,w)\,.
\end{eqnarray*}
\end{Proof}

\begin{Proposition} \label{PRO_braided_commutator_and_kernel}
Let $t\,\in\,X$ and $v\,\in\,\ker\partial_t$ be arbitrary. Then
$[e_t,\,v]_c\,\in\,\ker\partial_t$.
\end{Proposition}
\begin{Proof}
With Proposition \ref{PRO_skew-derivation_and_g_commute} one finds
\begin{eqnarray*}
\partial_t([e_t,\,v]_c) &=& \partial_t(e_t\,v) \;-\; \partial_t\big((g_t\, v)e_t\big)
\;\;=\;\; g_t\, v \;-\; g_t\,v \;\;=\;\; 0\,.
\end{eqnarray*}
\end{Proof}


\section{The Shift Group of a Nichols Algebra} \label{SEC_shifts}

Milinski and Schneider showed in \cite{MiS}, Theorem 5.8, that the grading of
a Nichols algebra is balanced, if $\Inn X$ is a Coxeter group and a certain type of
cocycle is given. In general, the grading of a Nichols algebra need not be
balanced, as we will see in the case of the 72-dimensional Nichols algebra
in subsection \ref{SEC_nichols_4B}.

Let $\B$ be a finite-dimensional Nichols-Algebra over the rack $X$ and cocycle $\chi$
with $\dim V_0\,=\,1$ (see Definition \ref{DEF_nichols_algebra}).
Recall that $\B$ is generated as an algebra by the elements $e_t\,\in\,\B(1)$,
$t\in X$.

It is a well-known fact that for each finite-dimensional Nichols algebra $\B$
over the quandle $X$ and the field $\K$ and for each $t\in X$ holds:

\begin{Lemma} \label{LEM_right_split}
Let $t\,\in\,X$ be arbitrary and $m\,=\,\nord e_t$. Then holds
$$ \B \;\cong\; (\ker \partial_t)\,\otimes\,(\K[e_t]/e_t^m)
\quad \textnormal{and} \quad
\B \;\cong\; (\K[e_t]/e_t^m)\,\otimes\,(\ker \partial^\op_t)\,. $$
\end{Lemma}
\begin{Proof}
1) This directly follows from Gra\~na's Freeness Theorem. 
We reproduce a short proof to compare it to part (2).

Let $s\in X\setminus\{t\}$ be arbitrary. Then by definition of the braided commutator
holds:
$$ e_t\,e_s \;=\; [e_t,\,e_s]_c \,+\, (g_t\,e_s)e_t\,.$$
For $[e_t,\,e_s]_c$ we use Proposition \ref{PRO_braided_commutator_and_kernel} to
see that $[e_t,\,e_s]_c\,\in\,\ker\partial_t$.
On the other hand, we have $\partial_t(g_t\, e_s)\,=\,0$, because $t\trid s\neq t$
for $s\neq t$ and any quandle $X$.
By induction we find that for each $v\in\B$ there are $v_j\in\ker \partial_t$
with
\beq
v \;=\; \sum_{j=0}^{m-1} v_j e_t^j\,.
\eeq
We now show that these $v_j$ are uniquely determined
(this is analog to Lemma 2.5 in \cite{MiS}):
Assume $\sum_{j=0}^{m-1}\,v_j e_t^j\,=\,0$. Apply $\partial_t$ $(m-1)$-times to
find $v_{m-1}\,=\,0$. Then apply $\partial_t$ $(m-2)$-times to see
$v_{m-2}\,=\,0$, induction.

2) Let $v\,\in\,\B$ be homogeneous with $\partial^\op_t(v)\,\neq\,0$.
By induction over the length, we can restrict to $v\,=\,u\,e_s$ with
$\partial^\op_t(u)\,=\,0$ but $\partial^\op_t(u\,e_s)\,\neq\,0$.
Set $w\,:=\,\chi(u,t)^{-1}\,u$. Then $w\,\in\,\ker\partial^\op_t$ and
\begin{eqnarray*}
\partial_t^\op(\underbrace{u\,e_s}_{=\,v}\,-\,e_t\,w) &=&
\chi(u,t)^{-1}\,u\,\underbrace{\partial^\op_{u\trid^{-1} t}\,e_t}_{=\,1}
\,-\,w\,-\,q_t^{-1}\,e_t\,\partial^\op_t\,w \;\;=\;\; 0\,,
\end{eqnarray*}
hence $v\,=\,e_t\,w\,-\,(v\,-\,e_t\,w)\,\in\,(\K[e_t]/e_t^m)\,\otimes\,(\ker \partial^\op_t)$.
Like in (1), linear independence is shown by applying $\partial_t^\op$.
\end{Proof}

One would expect that for each element $v\,\in\,\B(1)$, there is a decomposition
$\B\,=\,U\,\otimes\,\K[v]/v^{\nord v}$ similar to the one of Lemma \ref{LEM_right_split}.
This, however, is wrong: Take $v\,=\,e_1\,+\,e_2\,\in\,\na{3A}$. 
If $K$ is of characteristic $\neq\,2$, $v$ has nilpotency order $4$. If $\na{3A}$ would decompose
into a tensor product with factor $\K[v]/v^4$, its Hilbert series would be divisible
by $(4)_t$, which is not the case.

\begin{Proposition} \label{PRO_grades_of_coefficients}
Let $v\,\in\,\B(g)$ for some $g\,\in\,\Env X$ and $t\,\in\,X$ be arbitrary.
If we decompose $v$ into the
sum $\sum_{j=0}^{m-1}\,v_j\,e_t^j$ with Lemma \ref{LEM_right_split}, we have
$v_j\,\in\,\B(g\,t^{-j})$ for each $j$.
If decomposed into $v\,=\,\sum_{j=0}^{m-1}\,e_t^j\,v_j$, each $v_j\,\in\,\B(t^{-j}\,g)$.
\end{Proposition}
\begin{Proof}
Each summand $v_j e_t^j$ is itself $\Env X$-homogeneous (otherwise there would be
a non-trivial linear dependency). Due to uniqueness, all $v_j\,e_t^j$
are linearly independent, hence each $v_j\, e_t^j$ must be element
of $\B(g)\ni v$. Then $v_j\in \B(h)$ with $h=g\,t^{-j}$. The second
statement follows the same way.
\end{Proof}

\begin{Proposition} \label{PRO_shifts_are_isomorphisms}
Let $t,\,s\,\in\,X$ be arbitrary, $t\neq s$. The shift
\beq
\phi_t:\, \B\;\rightarrow\;\B,&& v\,=\,\sum_{j=0}^{m-1} v_j\,e_t^j
\;\mapsto\; v_{m-1}\,+\,\sum_{j=1}^{m-1} v_{j-1}\, e_t^j \\
&& \forall j:\, v_j\,\in\,\ker \partial_t
\eeq
is a well-defined linear isomorphism.
We call the group generated by the shifts $\phi_t$ the {\em shift group}
$\Phi(\B)$ (or just $\Phi$). By definition, its operation
on $\B$ is free.
\end{Proposition}
\begin{Proof}
$\phi_t$ is bijective because $\phi_t^m\,=\,\id_{\B}$ and $\phi_t$ obviously is linear.
\end{Proof}

\begin{Lemma} \label{LEM_phi_one_spans_B}
The $\Phi$-orbit of $1$ linearly spans $\B$.
\end{Lemma}
\begin{Proof}
By induction over the $\N_0$-degree $d$ of $v\,\in\,\B$.
For $d\,\in\,\{0,1\}$, this is clear. Assume $\Phi(1)$ spans the whole of
$\B(d)$. Let $w\,\in\,\B(d)$ and $t\,\in\,X$ be arbitrary.
By Lemma \ref{LEM_right_split}, $w$ decomposes into
$$ w \;\;=\;\; \sum_{j=0}^{m-1}\, w_j\, e_t^j $$
with $w_j\,\in\,\ker\partial_t$, $j\,=\,0\ldots (m-1)$.
Due to the grading, each $w_j$ can be chosen to be of length $d-j$.
Now we see
$$ w\,e_t \;\;=\;\; \sum_{j=0}^{m-2}\, w_j\, e_t^{j+1}
\;\;=\;\; \phi_t\left(\sum_{j=0}^{m-2}\, w_j\, e_t^j\right)\,. $$
$\sum_{j=0}^{m-2}\, w_j\, e_t^j$ is of length $d$ and can be spanned by
$\Phi(1)$. Hence, $w\,e_t$ can be spanned by $\Phi(1)$ as well.
\end{Proof}

As an (unused) corollary, we see that if $\Phi$ is finite,
$\B$ must be finite-dimensional. However, the converse is not true:
$\Phi(\na{4B})$ contains the infinite group $C_2\ast C_2$ as a subgroup
(see subsection \ref{SEC_nichols_4B}).

The dimension of the matrix algebra $\Alg \Phi$ spanned by the maps $\phi_t$
is bounded from above by $(\dim \B)^2$; and if $\B$ is infinite, it must
be infinite-dimensional as well, due to Lemma \ref{LEM_phi_one_spans_B}.
In the case of the Nichols algebra $\na{3A}$, we find that the dimension of
this shift algebra is $12$, which equals the dimension of $\B$. But we have
$\phi_t^2\,=\,\id$ in this algebra, so it cannot be $\N_0$-graded.
Hence, the shift algebra and $\B$ are not isomorphic.

In the case of $\na{4A}$, even the dimensions differ: $\B$ has dimension $36$,
but the shift algebra of $\B$ has dimension $648\,=\,18\cdot 36$.
For the 72-dimensional algebra $\na{4B}$ one finds that the shift algebra
has dimension $2592\,=\,36\cdot 72$.

\subsection{Case $n$ divides $m$} \label{SEC_n_divides_m}

In the following, let $\B$ be a Nichols algebra with $n\,|\,m$, $\Phi$ its
shift group and $\K\Phi$ the group algebra of $\Phi$. (Instead of $\K\Phi$
we might just as well take the shift algebra of $\B$.)

The evaluation at $1\,\in\,\B$ yields a linear map $\ev_1:\,\K\Phi\,\rightarrow\,\B$,
which by Lemma \ref{LEM_phi_one_spans_B} is surjective. In particular, we may
define subspaces $\K\Phi_g \,:=\, \ev_1^{-1}(\B(g))$.

\begin{Proposition} \label{PRO_phi_homogeneous}
If $\phi_t$ is restricted to $\B(g)$ for any $g\in \Inn X$, its image restricts
to $\B(gt)$ and hence yields a linear isomorphism $\B(g)\,\cong\,\B(gt)$.
\end{Proposition}
\begin{Proof}
We use the notations from Proposition \ref{PRO_shifts_are_isomorphisms}.
Let $v\in\B(g)$. Then $v_j\in\B(gt^{-j})$ and hence $v_{j-1} \,e_t^j\,\in\,
\B(gt)$ for each $j=1\ldots(m-1)$. Because of $n\,|\,m$ we have $t^m=1$ in $\Inn X$
and $v_{m-1}\in\B(h)$ with $h=gt^{-(m-1)}=gt$, so the sum $\phi_t(v)$ is in $\B(gt)$.
Apply $\phi_t$ $m$-times to see that $\phi_t|_{\B(g)}: \,\B(g)\rightarrow \B(gt)$ is
surjective and hence an isomorphism.
\end{Proof}

\begin{Corollary} \label{COR_nichols_balanced_for_n_divides_m}
The grading of $\B$ is balanced. The order of $\Inn X$ divides $\dim \B$.
\end{Corollary}
\begin{Proof}
$X$ generates $\Inn X$.
\end{Proof}

In particular, we see that $\K\Phi_g\cdot \K\Phi_h \,\subseteq\,\K\Phi_{gh}$, so
$\K\Phi/\ker \ev_1$ has a $G$-grading---it actually is just the $G$-grading
of $\B$ pushed to $\K\Phi/\ker \ev_1 \,\cong\,\B$.

\begin{Lemma} \label{LEM_G_is_quotient_of_Phi}
$\Inn X$ is a quotient of $\Phi$.
\end{Lemma}
\begin{Proof}
By Proposition \ref{PRO_phi_homogeneous}, the map $\gamma\,:=\,\deg\,\circ\,\ev_1:\,\Phi\,\rightarrow\,G$
is a surjective homomorphism.
\end{Proof}

Given a subset $X'$ of $X$, let $\B'$ be the subalgebra-with-one
of $\B$ generated by $\K X'\,\subseteq\,\B(1)$. If $X'$ is a subrack of $X$,
then $\B'$ is its corresponding Nichols subalgebra.

\begin{Proposition} \label{PRO_neutral_in_subalgebra}
Let $X'$ be a subrack of $X$ and $v\,\in\,\B'$.

1) If $\partial_t\,v\,=\,0$ for all $t\,\in\,X'$, then $v\,\in\,\B(0)$.

2) If $\partial_t^\op\,v\,=\,0$ for all $t\,\in\,X'$, then $v\,\in\,\B(0)$.
\end{Proposition}
\begin{Proof}
From $v\,\in\,\B'$ we know $\partial_s v\,=\,0\,=\,\partial_s^\op v$ for all $s\,\in\,X\setminus X'$,
so $v$ must be a multiple of $1$. 
Note that (1) actually holds for arbitrary subsets $X'$ of $X$ when $\B'$ is
defined as the sub-algebra generated by all $e_t$ with $t\,\in\,X'$; this is not
true for (2).
\end{Proof}

\begin{Lemma} \label{LEM_applying_granas_freeness_theorem}
Let $X'$ be a non-empty proper subrack of $X$ and $G'$ the subgroup of $\Inn X$
spanned by the operations $g_t:\,X\,\rightarrow\,X$ for all $t\,\in\,X'$.
Let $V'$ be the linear span of the elements $e_t$ with $t\,\in\,X'$ and
$U$ the linear span of all $e_s$ with $s\,\in\,X\setminus X'$. Then holds:
\begin{enumerate}
\item
$G'$ is the smallest subgroup of $G\,=\,\Inn X$ with $\delta(V')\,\subseteq\,\K G'\,\otimes\,V'$,
\item
$V'$ is $G'$-stable,
\item
$V'$ and $U$ are $\K G$-subcomodules, and
\item
$U$ is a $\K G'$-submodule.
\end{enumerate}
In particular, Gra\~na's Freeness Theorem (Theorem \ref{THE_grana_freeness}) applies and
we find
$$ \B\;\;\cong\;\;\left(\bigcap_{t\in X'}\,\ker\partial_t\right)\,\otimes\,\B'\,,$$
where $\B'$ is the Nichols sub-algebra generated by $X'$.
\end{Lemma}
\begin{Proof}
(1) holds by definition of $G'$. (2) holds because $X'$ is a subrack.
(3) is due to the diagonal comodule structure $\delta(e_t)\,=\,g_t\,\otimes\, e_t$
for all $t\in X$. To show (4), let $t\,\in\,X'$ and $s\,\in\,X\setminus X'$ be arbitrary.
Then $g_t(e_s)$ is a multiple of $e_{t\trid s}$. The element $t\trid s$
cannot be in $X'$ (otherwise $s$ would be in $X'$), so $g_t(e_s)\,\in\,U$.
\end{Proof}

Freeness theorems as Grana's allow for recursion, such that $\B$ can be
written as tensor product of terms of the form $\bigcap_{t\in X'}\,\ker\partial_t$
for ever decreasing subracks $X'$. Such a factorization induces a factorization of the
Hilbert series as well. In the case of the related Fomin-Kirillov algebras,
the analogous factorization has been conjectured in \cite{K_quadratic},
Conjecture 8.6, and has been proven by Fomin and Procesi in \cite{FP_quadratic}.
Their factors are subalgebras generated by transpositions $(i,\,n)$ for
fixed $n$ and $1\,\leq\,i\,<\,n$. This corresponds to the intersection
$\bigcap_{t\in X'}\,\ker\partial_t$ if $X$ is the rack generated by the transpositions
in $\mathfrak{S}_n$ and $X'$ the subrack generated by the transpositions of
$\mathfrak{S}_{n-1}\,<\,\mathfrak{S}_n$. The subalgebra generated by $X\setminus
X'$ is a subspace of $\bigcap_{t\in X'}\,\ker\partial_t$, but not
necessarily all of it.

\begin{Proposition} \label{PRO_kernel_homogeneous_basis}
Let $X'$ be a non-empty subset of $X$. Let $G\,=\,\Env X$ or any quotient of $\Env X$.
Then $\bigcap_{t\in X'}\,\ker\partial_t\,\subseteq\,\B$ has a $G$-homogeneous basis.
\end{Proposition}
\begin{Proof}
Each $\partial_t$ is a $G$-graded map, hence $\ker\partial_t$ is a $G$-graded
sub-$\B$-module of $\B$; same for their intersection
$\bigcap_{t\in X'}\,\ker\partial_t$. Each graded submodule has a homogeneous
basis, see e.g.\ section 2.1 in \cite{NO_graded_rings}.
\end{Proof}

\section{The Modified Shift Groups} \label{SEC_modified_shifts}

Each shift $\phi_t$ can be written in the form
$$ \phi_t(v) \;\;=\;\; \frac{1}{(1+q_t)(1+q_t+q_t^2)\cdots(1+q_t+\cdots+q_t^{m-2})} \cdot \partial_t^{m-1} v \;+\; v\,e_t $$
by inserting the decomposition of $v$ into the right-hand side.
We modify this definition by removing the leading factor and subtracting the
braided commutator for one variant; and by issuing the opposite braided derivative
for another.

\begin{Definition} \label{DEF_modified_shift}
Let $\B$ be a Nichols algebra with rack $X$ and $t\,\in\,X$ arbitrary.
Define the {\em modified shifts}
$$ \psi_t(v) \;:=\; \partial_t^{m-1}(v) \;+\; e_t\,g_t^{-1}(v) $$
and
$$ \xi_t(v) \;:=\; (\partial^\op_t)^{m-1}(v) \;+\; e_t\,v $$
and the corresponding {\em modified shift groups} $\Psi$ and $\Xi$,
generated by all $\psi_t$ (respectively $\xi_t$) with $t\,\in\,X$.
If $X'$ is a subset of $X$, define $\Psi|_{X'}$ and $\Xi|_{X'}$
to be the groups generated by all $\psi_t$ (respectively $\xi_t$)
with $t\,\in\,X'$.
\end{Definition}

\begin{Proposition} \label{PRO_modified_shift}
Let $s,\,t\,\in\,X$ and $g\,\in\,G$ be arbitrary, where $G$ is any quotient of $\Env X$.
Then holds:
\begin{enumerate}
\item $\psi_t$ is a linear isomorphism.
\item If $t^m\,=\,e$ in $G$, then $\psi_t$ maps $\B(g)$ to $\B(gt)$.
\item If $t\,\neq\,s$, then $\psi_t(\ker\partial^\op_s)\,=\,\ker\partial^\op_s$.
\item The $\Psi$-orbit of $1$ linearly spans $\B$.
\end{enumerate}
\end{Proposition}
\begin{Proof}
{\bf 1.} Linearity is obvious. Now assume $v\,\in\,\ker\psi_t\setminus\{0\}$. Then
$\partial_t^{m-1}(v)\,=\,-e_t\,g_t^{-1}(v)$. Use Lemma \ref{LEM_right_split} to decompose
$v\,=\,\sum_{j=0}^{m-1}\,v_j\,e_t^j$ with $v_j\,\in\,\ker\partial_t$. Inserting this
yields $v_{m-1}\,=\,\sum\,\lambda_j\,e_t\,g_t^{-1}(v_j)\,e_t^j$ for some
$\lambda_j\,\in\,\K\setminus\{0\}$.
Using the braided commutator, we find
$$ e_t\,g_t^{-1}(v_j)\,e_t^j \;\;=\;\; v_j\,e_t^{j+1}
\;+\; \underbrace{[e_t,\,g_t^{-1}(v_j)]_c}_{\in\,\ker\partial_t}\,e_t^j $$
by Propositions \ref{PRO_skew-derivation_and_g_commute} and \ref{PRO_braided_commutator_and_kernel}.
We know $v_{m-1}\,\in\,\ker\partial_t$, and by comparing the coefficients of $e_t^j$,
we find:
\begin{eqnarray*}
v_{m-1} &=& \lambda_0\,[e_t,\, g_t^{-1}(v_0)]_c \\
v_{j-1} &=& -\,[e_t,\,g_t^{-1}(v_j)]_c \qquad \textnormal{if $1\,\leq\,j\,\leq\,m-1$} 
\end{eqnarray*}
In particular, the minimal length of the $\N_0$-homogeneous components of $v_{j-1}$
is at least one plus the minimal length of $v_j$, and the minimal length of $v_{m-1}$ is
at least one plus the minimal length of $v_0$. This cannot be, hence all $v_j$ are zero.

{\bf 2.} Let $v\,\in\,\B(g)\setminus\{0\}$ be arbitrary. Then the degree of
$\partial_t^{m-1}(v)\,=\,gt^{1-m}\,=\,gt$ and the degree of $e_t\,g_t^{-1}(v)$
also is $t\,t^{-1}\,gt\,=\,gt$. Thereby, $\psi_t(v)$ is homogeneous of degree $gt$.

{\bf 3.} Let $v\,\in\,\ker\partial_s^\op$ be arbitrary. From Proposition
\ref{PRO_left_derivative} we conclude that $\partial_t^{m-1}(v)\,\in\,\ker\partial_s^\op$.
For the other summand, we find
$$ \partial_s^\op(e_t\,g_t^{-1}(v)) \;\;=\;\;
q\,e_t\,\partial_{t\,\trid^{-1}\,s}^\op\,g_t^{-1}(v)\,. $$
$v$ is in $\ker\partial_s^\op$, so $g_t^{-1}(v)$ is in $\ker\partial_{t\,\trid^{-1}\,s}^\op$.

{\bf 4.} We may use the same induction as in the proof of Lemma \ref{LEM_phi_one_spans_B}:
We write
$$ e_t \,v \;=\; \psi_t(g_t\,v) \,-\, \partial_t^{m-1}(g_t\,v) $$
and note that both summands on the right hand side are in the linear span of $\Phi(1)$.
\end{Proof}

Property (3) above is something not found in
the shifts introduced in Proposition \ref{PRO_shifts_are_isomorphisms}:
In the Nichols algebra $\na{3A}$ choose
$v\,=\,e_3$. Then $v\,\in\,\ker\partial_1^\op$, but $\phi_2(v)\,=\,e_3\,e_2\,
\notin\ker\partial_1^\op$, while
$\psi_2(v)\,=\,-e_2\,e_1\,\in\ker\partial_1^\op$. This justifies
the introduction of the modified shifts.

\begin{Proposition} \label{PRO_modified_shift_xi}
Let $s,\,t\,\in\,X$ and $g\,\in\,G$ be arbitrary, where $G$ is any quotient of $\Env X$.
Then holds:
\begin{enumerate}
\item $\xi_s$ is a linear isomorphism.
\item If $s^m\,=\,e$ in $G$, then $\xi_s$ maps $\B(g)$ to $\B(s\,g)$.
\item If $t\,\neq\,s$, then $\xi_s(\ker\partial_t)\,=\,\ker\partial_t$.
\item The $\Xi$-orbit of $1$ linearly spans $\B$.
\end{enumerate}
\end{Proposition}
\begin{Proof}
The proof is very similar to the proof of Proposition \ref{PRO_modified_shift},
but differs in details.

{\bf 1.} Again, linearity is obvious. Let $v\,\in\,\ker\xi_s\setminus\{0\}$
be arbitrary and use the right-hand decomposition in Lemma \ref{LEM_right_split},
$v\,=\,\sum_{j=0}^{m-1}\,e_s^j\,v_j$ with $v_j\,\in\,\ker\partial^\op_s$. Inserting this
into $\xi_s(v)\,=\,0$ yields
$$ v_{m-1}\,=\,\sum_{j=0}^{m-2}\,\mu_j\,e_s^{j+1}\,v_j $$
for some $\mu_j\,\in\,\K\setminus\{0\}$. Due to the linear independence in the decomposition
in Lemma \ref{LEM_right_split}, we conclude $v_j\,=\,0$ for all $j$ and thus $v\,=\,0$.

{\bf 2.} Straightforward, see Proposition \ref{PRO_modified_shift}.

{\bf 3.} Let $v\,\in\,\ker\partial_t$ be arbitrary. Due to Proposition
\ref{PRO_left_derivative} $(\partial_s^\op)^{m-1}$ and $\partial_t$ commute,
so $(\partial^\op_s)^{m-1}(v)\,\in\,\ker\partial_t$.
The other summand vanishes due to $\partial_t(e_s\,v) \,=\, e_s\,\partial_t(v) \,=\, 0$.

{\bf 4.} Analog to Proposition \ref{PRO_modified_shift}.
\end{Proof}

We will mainly use the shifts $\xi_t$ in the following, due to the special form
of Gra\~nas Freeness Theorem we are using.

\begin{Lemma} \label{LEM_kernel_balanced}
Let $X'$ be a non-empty subset of $X$. Let $G$ be a quotient of $\Env X$
with $s^m\,=\,e$ for all $s\,\in\,X\setminus X'$,
such that $X\setminus X'$ still generates $G$.
Then $\bigcap_{t\in X'}\,\ker\partial_t$ is $G$-balanced.
\end{Lemma}
\begin{Proof}
By Proposition \ref{PRO_kernel_homogeneous_basis}, $K\,:=\,\bigcap_{t\in X'}\,\ker\partial_t$
has a $G$-homogeneous basis. For any $s\,\in\,X\setminus X'$ and $g\,\in\,G$,
$\xi_s$ will map $K\,\cap\,\B(g)$ to $K\,\cap\,\B(s\,g)$ by Proposition
\ref{PRO_modified_shift_xi}, so $\dim (K\,\cap\,\B(g))\,=\,\dim(K\,\cap\,\B(s\,g))$.
The Lemma now follows from the assumption that $G$ is generated by $X\setminus X'$.
\end{Proof}

\begin{Theorem} \label{THE_dimension_divisibility}
Assume $n\,|\,m$. Let $X'$ be a non-empty proper subrack of $X$, and assume
that $X\setminus X'$ still generates $\Inn X$. Set $\B'$ to be the
sub-Nichols algebra generated by $X'$. Then holds
\begin{equation} \label{EQU_dimension_divisibility}
\dim\B \;\;=\;\; \#\Inn X\cdot \dim\B'\cdot \dim\left(\B(e)\,\cap\,\bigcap_{t\in X'}\,\ker\partial_t\right)\,.
\end{equation}
\end{Theorem}
\begin{Proof}
The proof follows directly from Gra\~nas Freeness Theorem
by applying Lemma \ref{LEM_kernel_balanced}
to Lemma \ref{LEM_applying_granas_freeness_theorem}
(note that $g_t^m\,=\,e$ holds for the inner group if and only if $n\,|\,m$).
\end{Proof}

Assume $X$ is indecomposable and consists of at least three elements
(there is no indecomposable quandle of two elements).
Choose $X'\,=\,\{t\}\,\subsetneq\,X$, thus $\B'\,\cong\,\K[t]/t^n$.
Due to irreducibility, there must exist $r,\,s\,\in\,X\setminus X'$
with $r\,\trid\,s\,=\,t$, so $g_t\,=\,g_r\,g_s\,g_r^{-1}$ and $\Inn X$
is generated by $X\setminus X'$. Applying Theorem
\ref{THE_dimension_divisibility}, we find $\#\Inn X\,\cdot n\,|\,\dim \B$.

\begin{Example} \label{EXA_nichols_transpositions_orders}
From the examples of subsection \ref{SEC_tables}, four 
Nichols algebras fulfill $n\,|\,m$:
$$
\!\!\!\!
\begin{array}{|r|c|r|r|c|r|r|} \hline
  \B & m & \dim \B & \# G  & \B' & \# G\cdot \dim \B' &
\!\!\!\! \frac{\dim \B}{\#G\cdot \dim \B'} \!\!\!\! \\ \hline
  \na{3A} & 2 &        12 &   6 & \K[t]/t^2 &     12 &   1 \\ 
  \na{4C} & 3 &     5,184 &  12 & \K[t]/t^3 &     36 & 144 \\
  \na{6A} & 2 &       576 &  24 &   \na{3A} &    288 &   2 \\
 \na{10A} & 2 & 8,294,400 & 120 &   \na{6A} & 69,120 & 120 \\ \hline
\end{array}
$$
where we use $G\,=\,\Inn X$.

It is still unclear, whether $\na{15A}$ (the Nichols algebra of the
transpositions in the symmetric group $\mathfrak{S}_6$ with constant cocycle $-1$)
is finite-dimensional or not.
If it is, its dimension must be divisible by
$$ \# \mathfrak{S}_6 \cdot \dim\na{10A} \;=\; 720\cdot 8,294,400 \;=\; 5,971,968,000\,. $$
Taking a look at the quotients $\frac{\dim \B}{\#G\cdot \dim \B'}$ in the above table,
one might guess that $\dim \na{15A}$ will probably be about at least another factor
$720$ larger, and thus divisible by $4,299,816,960,000$.
\end{Example}


\section{General Case} \label{SEC_general_case}

We remember from Lemma \ref{LEM_G_is_quotient_of_Phi} that $\Inn X$ is a
quotient of the shift group $\Phi$ if $n\,|\,m$. Moreover, this induces a
$G$-grading on $\K\Phi$ (which is balanced if $\K\Phi$ is finite-dimensional),
which in turn induces a balanced $G$-grading on $\B$. One might ask, how
to generalize this idea to the case $n\,\nmid\,m$.

Let $\ev_1:\,\K\Phi\,\rightarrow\,\B$
be the evaluation at $1\,\in\,\B$. From Lemma
\ref{LEM_phi_one_spans_B} we know that $\ev_1$ is a surjective linear
map. $\ev_1$ is neither an algebra homomorphism, nor is its kernel
an ideal of $\K\Phi$; still, there is an identification of $\K\Phi/\ker\ev_1$
and $\B$ as linear spaces. Assume there is a surjective homomorphism
$\pi:\,\Phi\,\rightarrow\,G$ to some finite quotient $G$ of $\Env X$.
Define $\Phi_g\,:=\,\pi^{-1}(g)$ and
$U_g\,:=\,\K\Phi_g\,\subseteq\,\K\Phi$, such that $\K\Phi\,=\,
\oplus_{g\in G}\,U_g$. Choose a system of representatives
$\phi_g\,\in\,U_g\setminus\{0\}$ and define the translations
$\tau_g:\,\K\Phi\,\rightarrow\,\K\Phi$, $\phi\,\mapsto\,\phi_g\,\phi$.
Each $\tau_g$ is a linear isomorphism of $\K\Phi$ and $\tau_g(U_h)\,=\,U_{gh}$
for each $g,\,h\,\in\,G$. Now assume $\phi\,\in\,\ker\ev_1$. Then
$\tau_g(\phi)(1)\,=\,\phi_g(\phi(1))\,=\,0$, hence
$\tau_g(\ker\ev_1)\,=\,\ker\ev_1$.
We may therefore define linear maps $\tilde\tau_g:\,\K\Phi/\ker\ev_1\,
\rightarrow\,\K\Phi/\ker\ev_1$ with $\tilde\tau_g(U_h/\ker\ev_1)\,=\,
U_{gh}/\ker\ev_1$.
Obviously, we have $\K\Phi/\ker\ev_1\,=\,\sum_{g\in G}\, U_g/\ker\ev_1$.
Assume this sum is direct. Then the isomorphisms $\tilde\tau_g$ show
that this grading is balanced, and hence $\#G$ divides $\dim \B$.
So there currently are two open questions to transfer the results
of Subsection \ref{SEC_n_divides_m} to the general case:
\begin{enumerate}
\item Is $G$ a quotient of $\Phi$?
\item Is the sum $\sum_{g\in G}\, U_g/\ker\ev_1$ direct?
\end{enumerate}
We will now concentrate on $G\,=\,C_k$, which is a quotient of $\Env X$
by $t\,\mapsto\,[1]_k\,\in\,C_k$ for all $t\,\in\,X$.

\subsection{Factors in the Hilbert Series} \label{SEC_hilbert}

Each Nichols Algebra $\B$ is $\Z$-graded. Taking quotients, we find
$C_k$-gradings of $\B$ for each $k\,>\,1$.

\begin{Lemma} \label{LEM_Cm_balanced}
$\B$ is $C_m$-balanced. In particular, we have for each $[k]_m\,\in\,C_m$
$$ \sum_{\begin{array}{c} j\in\N_0,\\ j \equiv k\; (\mathrm{mod}\, m)\end{array}}
\!\!\!\!\!\!\!\!\!\!\!\! \dim \B(j) \;\;=\;\; \frac{1}{m}\,\dim \B\,. $$
\end{Lemma}
\begin{Proof}
Let $t\,\in\,X$ be arbitrary. We may define $\phi_t$ as in Proposition
\ref{PRO_shifts_are_isomorphisms} (or, equivalently, any of the shifts of
Definition \ref{DEF_modified_shift})
and find that $\phi_t$ maps $\B(j)$ to $\B(j+1)\,\oplus\,\B(j-m+1)$.
Hence, $\phi_t$ is a linear isomorphism between $\B([j]_m)$ and
$\B([j+1]_m)$.
\end{Proof}

Clearly, if $j\,|\,k$ and $\B$ is $C_k$-balanced, then $\B$ is $C_j$-balanced
as well. The following table shows for some Nichols algebras $\B$ those
$k\,>\,1$ such that $\B$ is $C_k$-balanced.

\begin{center}
\begin{tabular}{|r|c|c|r|l|} \hline
 \;  $\B\phantom{/B/C}$ \; & $n$ & $m$ & $\dim \B$ & $C_k$-balanced for $\ldots$ \\ \hline
 \;  $\K[t]/t^m$ \; & $1$     & $m$ &       $m$ & \; $k\;|\;m$ \\
 \;  $\na{3A}$ \; & $2$     & $2$ &        12 & \; $k\,=\,$ 2, 3\phantom{, 4, 5, 6, 7} \; \\ 
 \;  $\na{3B}$ \; & $2$     & $3$ &       432 & \; $k\,=\,$ 2, 3, 4, \phantom{5, }6\phantom{, 7} \\ 
 \;  $\na{4A}$ \; & $3$     & $2$ &        36 & \; $k\,=\,$ 2, 3\phantom{, 4, 5, 6, 7} \\
 \;  $\na{4B}$ \; & $3$     & $2$ &        72 & \; $k\,=\,$ 2, 3, \phantom{4, 5, }6\phantom{, 7} \\
 \;  $\na{4C}$ \; & $3$     & $3$ &     5,184 & \; $k\,=\,$ 2, 3, 4, \phantom{5, }6\phantom{, 7} \\
 \;  $\na{5*}$ \; & $4$     & $2$ &     1,280 & \; $k\,=\,$ 2, \phantom{3, }4, 5\phantom{, 6, 7} \\
 \;  $\na{6*}$ \;             & $2/2/4$ & $2$ &       576 & \; $k\,=\,$ 2, 3, 4\phantom{, 5, 6, 7} \\
 \;  $\na{7*}$ \; & $6$     & $2$ &   326,592 & \; $k\,=\,$ 2, 3, \phantom{4, 5, }6, 7 \\
 \; $\na{10*}$ \; & $2$     & $2$ & 8,294,400 & \; $k\,=\,$ 2, 3, 4, 5, 6\phantom{, 7} \\ \hline
\end{tabular}
\end{center}

\begin{Lemma} \label{LEM_Ck_balanced_and_Hilbert}
A finite-dimensional Nichols algebra $\B$ is $C_k$-balanced 
(as a quotient of its $\Z$-grading) if
and only if $(k)_t\,:=\,\sum_{j=0}^{k-1}\,t^j$ is a divisor of
the Hilbert series $\mathcal{H}_\B(t)$ of $\B$, such that the quotient
polynomial has integer coefficients only.
\end{Lemma}
\begin{Proof}
To simplify notation, if $p$ is a polynomial, set $p_j$ to be the
coefficient of $t^j$ in $p(t)$ (or zero if $j\,<\,0$)
and $b_j\,:=\,(\mathcal{H}_\B)_j\,=\,\dim\B(j)$ for any $j\,\in\,\N_0$.

``$\Rightarrow$'': Set $p_j\,:=\,0$ for each $j\,<\,0$ and inductively define
$p_j\,:=\,b_j\,-\,\sum_{i=1}^{k-1}\,p_{j-i}$, hence $b_j\,-\,b_{j-1}\,=\,p_j\,-\,p_{j-k}$.
Let $d\,\in\,\N_0$ be such that $d\cdot k$ is larger than the top degree of $\B$.
Then summation of the previous equation yields for each $0\,\leq\,l\,\leq k\,-\,1$
a telescoping sum
\begin{eqnarray*}
\sum_{0\,\leq\,j\,\leq\,d}b_{jk+l}\,-\,\sum_{0\,\leq\,j\,\leq\,d}b_{jk+l-1} &=&
\sum_{0\,\leq\,j\,\leq\,d}\big(p_{jk+l}\,-\,p_{jk+l-k}\big) \\
&=& -\,p_{l-k}\,+\,p_{dk+l}\,.
\end{eqnarray*}
$\B$ is $C_k$-balanced by assumption, so the two sums on the left hand side must
sum to the same value (namely $\frac{1}{k}\dim\B$).
$p_{l-k}$ is zero by definition ($l-k<0$), hence $p_{dk+l}$
is zero as well. This shows that $p$ actually is a polynomial, and by definition its
coefficients are integers.
From $b_j\,=\,\sum_{i=0}^{k-1}\,p_{j-i}$,
we also see $\mathcal{H}_\B(t)\,=\,(k)_t \cdot p(t)$.

``$\Leftarrow$'': Let $\mathcal{H}_\B(t)\,=\,(k)_t\cdot p(t)$ for some polynomial $p$
with integer coefficients. We have $b_j\,=\,\sum_{i=0}^{k-1}\,p_{j-i}$
and therefore for each $[l]_k\,\in\,C_k$
\begin{eqnarray*}
\dim \B([l]_k) &=& \!\!\!\!\!\!\!\!\!\!\!\!
\sum_{\begin{array}{c} j\in\N_0,\\ j \equiv l\; (\mathrm{mod}\, k)\end{array}}
\!\!\!\!\!\!\!\!\!\!\!\! \dim \B(j) \;\;=\;\; \!\!\!\!\!\!\!\!\!\!\!\!
\sum_{\begin{array}{c} j\in\N_0,\\ j \equiv l\; (\mathrm{mod}\, k)\end{array}}
\!\!\!\!\!\!\!\!\!\!\!\! \sum_{i=0}^{k-1}\,p_{j-i}
\;\;=\;\;\sum_{j\in\N_0}\,p_j\,,
\end{eqnarray*}
which does not depend on $[l]_k$, so $\B$ is $C_k$-balanced.
\end{Proof}

From Lemmas \ref{LEM_Cm_balanced} and \ref{LEM_Ck_balanced_and_Hilbert} follows
that $(m)_t$ is a divisor of the Hilbert series of $\B$. This result is well-known
and can be seen directly from any of the Freeness Theorems applied to a trivial
subrack.

\begin{Theorem} \label{THE_hilbert_divisibility}
Let $\B$ be a finite-dimensional Nichols algebra over a rack $X$ and a 2-cocycle
of order $m$. Let $X'$ be a non-empty proper subrack of $X$ and $\B'$ its corresponding
Nichols sub-algebra of $\B$. Then the Hilbert series $\mathcal{H}_\B(t)$ is
divisible by $(m)_t\cdot \mathcal{H}_{\B'}(t)$.
\end{Theorem}
\begin{Proof}
We use the notation of Lemmas \ref{LEM_Cm_balanced} and \ref{LEM_Ck_balanced_and_Hilbert}.
Let $t\,\in\,X$ be arbitrary. We have seen in the proof of Lemma \ref{LEM_Cm_balanced},
that $\phi_t$ is a linear isomorphism between $\B([j]_m)$ and $\B([j+1]_m)$
for all $j$; so is $\xi_t$ from Proposition \ref{PRO_modified_shift_xi}.
Applying the same techniques of the proof of Theorem \ref{THE_dimension_divisibility}
to the quotient $C_m$-grading yields the proposition.
\end{Proof}

\begin{Corollary}
Let $\B$ be a finite-dimensional Nichols algebra over an indecomposable
rack $X$ with $\# X\,\geq\,3$ and a 2-cocycle of order $m$. Then its Hilbert
series $\mathcal{H}_\B(t)$ is divisible by $(m)_t^2$.
\end{Corollary}

\begin{Example}
We know that there is a sequence of embeddings of quandles
$$ \{t\}\;\;\hookrightarrow\;\;Q_{3,1}\;\;\hookrightarrow\;\;Q_{6,1}\;\;
\hookrightarrow\;\;Q_{10,1} $$
associated to the Nichols-algebra-embeddings
$$ \K[t]/t^2\;\;\hookrightarrow\;\;\na{3A}\;\;\hookrightarrow\;\;\na{6A}
\;\;\hookrightarrow\;\;\na{10A}\,. $$
In addition to this, one easily sees that $Q_{6,1}\setminus Q_{3,1}$
still generates $\Inn Q_{6,1}$ and $Q_{10,1}\setminus Q_{6,1}$ still
generates $\Inn Q_{10,1}$. Applying Theorem \ref{THE_hilbert_divisibility}
three times now shows that $(2)_t^4$ is a factor of $\mathcal{H}_{\na{10A}}(t)$.
Following Example \ref{EXA_nichols_transpositions_orders}, we conclude that
$(2)_t^5$ must be a factor of $\mathcal{H}_{\mathcal{B}(Q_{15,7},\,-1)}(t)$,
if this Nichols algebra is finite dimensional.
\end{Example}

\subsection{The 72-dimensional Nichols Algebra} \label{SEC_nichols_4B}

We now concentrate on one example with $n\,\nmid\,m$, the 72-dimensional
Nichols algebra $\na{4B}$ first introduced in \cite{G_low}.

Let $X=(\{1,2,3,4\},\,\trid)$ be the quandle with operation
\begin{center}
\begin{tabular}{c|cccc}
 $\trid$ & 1 & 2 & 3 & 4 \\ \hline
 1 & 1 & 3 & 4 & 2 \\
 2 & 4 & 2 & 1 & 3 \\
 3 & 2 & 4 & 3 & 1 \\
 4 & 3 & 1 & 2 & 4
\end{tabular}
\end{center}
Choose $\dim V_0\,=\,1$ and as cocycle choose the constant cocycle $\chi\,=\,-1$
over any field $\K$ of characteristic $\neq 2$.
The resulting Nichols Algebra $\B$ has dimension $72$
(\cite{G_low}), a possible basis is given by the following
products, written in syntax notation:
$$
\big[\;e_1\;\big]\;\big[\;e_2\;[\;e_1\;]\;\big]\;\big[\;e_3\,e_2\,e_1\;\big]\;
\big[\;e_3\;[\;e_2\;]\;\big]\;\big[\;e_4\;\big]
$$
(each argument in square brackets is optional). Its relations are generated
by the relations
\begin{eqnarray*}
0 &=& e_t^2 \qquad \forall t\,\in\, X \\
0 &=& e_r\,e_s\,+\,e_s\,e_t\,+\,e_t\,e_r \\
&& \qquad \qquad \qquad \forall (r,s,t)\,\in\,
\{(4,3,2),\,(4,2,1),\,(4,1,3),\,(3,1,2)\} \\
0 &=& (e_3\,e_2\,e_1)^2\,+\,(e_2\,e_1\,e_3)^2\,+\,(e_1\,e_3\,e_2)^2
\end{eqnarray*}
The inner group of $X$ is isomorphic to the alternating group $A_4$.
With respect to this grading, $\B$ is not balanced:
\begin{eqnarray*}
\begin{array}{|r||cccc|} \hline
g\in A_4  & ()  & (1,3)(2,4) & (1,4)(2,3) & (1,2)(3,4) \\ \hline
\dim \B(g) &  12  &     4      &     4      &     4     \\ \hline \hline
g\in A_4  & (1,2,3) & (1,3,4) & (1,4,2) & (2,3,4) \\ \hline
\dim \B(g) &    6    &    6    &    6    &    6    \\ \hline \hline
g\in A_4  & (1,3,2) & (1,4,3) & (1,2,4) & (2,4,3) \\ \hline
\dim \B(g) &    6    &    6    &    6    &    6    \\ \hline
\end{array}
\end{eqnarray*}
(Elements in cycle notation; calculations have been performed with Rig,
see \cite{Rig}.)
As one sees, the dimension is preserved
by conjugation; this is due to the operation of $\Env X$ on $\B$, which
conjugates the grading.

The grading of $\B$ with respect to the enveloping group $\Env X$ of $X$
must be unbalanced, because $\Env X$ is infinite. Indeed, the two elements
$g_3 g_2 g_1$ and $g_2 g_3 g_2 g_1 g_3 g_4\,\in\,\Env X$ fulfill $\dim \B(g) = 5$,
eight elements have $\dim \B(g) = 3$, another eight elements $2$, $22$ elements
have $\dim\B(g) = 1$ (including the identity element) and the remaining
elements $0$. One would therefore ask, whether there is a quotient $G$ of
$\Env X$, such that $\B$ is $G$-balanced and $G$ is large enough
to have $\Inn X$ as a quotient itself. However:

\begin{Proposition}
There is no quotient $G$ of $\Env X$, such that $\Inn X$ is a quotient of $G$ and
$\B$ is $G$-balanced.
\end{Proposition}
\begin{Proof}
Taking a quotient cannot lower the dimensions of the grade. For $g=g_3 g_2 g_1\in \Env X$
we therefore find, that the dimension of the grade of the image of $g$ under
the canonical projection $\Env X\rightarrow G$ (and hence for each element in $G$)
must be at least $5$. By hypothesis, $\Inn X$ is a quotient of $G$ and therefore
fulfills $\dim \B(g)\geq 5$ for each $g\in \Inn X$: Contradiction.
\end{Proof}

\begin{Proposition}
The shift group $\Phi$ of $\na{4B}$ is infinite in characteristic $0$.
\end{Proposition}
\begin{Proof}
The endomorphism $\phi_1\,\phi_2$ has a Jacobi normal form with eight blocks
of each of the three types
$$
\begin{pmatrix}
1 & 1 & 0 \\
0 & 1 & 1 \\
0 & 0 & 1
\end{pmatrix}
\qquad
\begin{pmatrix}
\rho & 1 & 0 \\
0 & \rho & 1 \\
0 & 0 & \rho
\end{pmatrix}
\qquad
\begin{pmatrix}
\bar\rho & 1 & 0 \\
0 & \bar\rho & 1 \\
0 & 0 & \bar\rho
\end{pmatrix}
$$
where $\rho$ and $\bar\rho$ are different third roots of unity.
From this decomposition one sees that $\phi_1\,\phi_2$ has infinite order
if $\K$ is of characteristic $0$, and thus
$\langle \phi_1,\,\phi_2\rangle\,\cong\,C_2\ast C_2$.
\end{Proof}

{\bf Acknowledgements.}
The author wants to thank Istv\'an Heckenberger and Leandro Vendramin
for their valuable hints and corrections.


\end{document}